\documentclass{amsart}
\usepackage{newlfont,amsfonts,amssymb,
amsmath,amsthm,amsgen,amscd,wasysym,comment}
\newcommand{\eeqref}[1]{equation \eqref{#1}}

\newcommand{\K}{\ensuremath{\mathbb{K}}}

\newcommand{\SO}{\mathrm{SO}}
\newcommand{\1}{\mathbb{I}}

\newcommand{\ups}{\Upsilon}
\newcommand{\e}{\mathrm{e}}

\newcommand{\bcase}{\begin{case}}
\newcommand{\ecase}{\end{case}}

\newcommand{\bclaim}{\begin{claim}}
\newcommand{\eclaim}{\end{claim}}

\newcommand{\bstep}{\begin{step}}
\newcommand{\estep}{\end{step}}

\newcommand{\bhlem}{\begin{hlem}}
\newcommand{\ehlem}{\end{hlem}}

\newcommand{\bleer}{\begin{leer}}
\newcommand{\eleer}{\end{leer}}
\newcommand{\bde}{\begin{de}}
\newcommand{\ede}{\end{de}}

\newcommand{\ol}{\overline}

\newcommand{\bs}{\begin{satz}}
\newcommand{\es}{\end{satz}}
\newcommand{\btheo}{\begin{theo}}
\newcommand{\etheo}{\end{theo}}
\newcommand{\bfolg}{\begin{folg}}
\newcommand{\efolg}{\end{folg}}
\newcommand{\blem}{\begin{lem}}
\newcommand{\elem}{\end{lem}}
\newcommand{\bnote}{\begin{note}}
\newcommand{\enote}{\end{note}}
\newcommand{\bprf}{\begin{proof}}
\newcommand{\eprf}{\end{proof}}
\newcommand{\bd}{\begin{displaymath}}
\newcommand{\ed}{\end{displaymath}}
\newcommand{\be}{\begin{eqnarray*}}
\newcommand{\ee}{\end{eqnarray*}}
\newcommand{\eeqa}{\end{eqnarray}}
\newcommand{\beqa}{\begin{eqnarray}}
\newcommand{\bi}{\begin{itemize}}
\newcommand{\ei}{\end{itemize}}
\newcommand{\bnum}{\begin{enumerate}}
\newcommand{\enum}{\end{enumerate}}

\newcommand{\eps}{\epsilon}

\newcommand{\beq}{\begin{equation}}
\newcommand{\eeq}{\end{equation}}
\newcommand{\einhalb}{\frac{1}{2}}
\newcommand{\rr}{\mathbb{R}}

\newcommand{\vf}{\varphi}
\newcommand{\earr}{\end{array}\]}
\newcommand{\barr}{\[\begin{array}}
\newcommand{\bvec}{\left(\begin{array}{c}}
\newcommand{\evec}{\end{array}\right)}

\newcommand{\sumi}{\sum_{i=1}^{n-2}}

\newcommand{\sumij}{\sum_{i,j=1}^{n-2}}

\newcommand{\+}{\oplus}

\newcommand{\laso}{\mathfrak{so}}

\newcommand{\w}{\omega}

\newcommand{\s}{\sigma}

\newcommand{\bbem}{\begin{bem}}
\newcommand{\ebem}{\end{bem}}
\newcommand{\bbez}{\begin{bez}}
\newcommand{\ebez}{\end{bez}}
\newcommand{\bbsp}{\begin{bsp}}
\newcommand{\ebsp}{\end{bsp}}

\newcommand{\trace}{\mathsf{tr }}

\newcommand{\ro}{\mathsf{P}}

\theoremstyle{definition}
\newtheorem{de}{Definition}
\newtheorem{bem}{Remark}
\newtheorem{bez}{Notation}
\newtheorem{bsp}{Example}
\theoremstyle{plain}
\newtheorem{lem}{Lemma}
\newtheorem{satz}{Proposition}
\newtheorem{folg}{Corollary}
\newtheorem{theo}{Theorem}

\begin{document}

\bibliographystyle{abbrv}


\title{Conformal pure radiation with parallel rays}
\author{Thomas Leistner}\address[Leistner]{School of Mathematical Sciences, University of Adelaide, SA 5005, Australia} \email{thomas.leistner@adelaide.edu.au}\thanks{This work was supported by the first author's Start-Up-Grant of the Faculty of
Engineering, Computer and Mathematical Sciences of the University of
Adelaide and by 
grants NN201 607540 and NN202 104838 of
the  Polish Ministry of Research and Higher Education.}
\author{Pawe\l~ Nurowski} \address[Nurowski]{Wydzia\l ~Fizyki, Instytut Fizyki Teoretycznej,
Uniwersytet Warszawski\\ ul. Ho\.za 69, 00-681 Warszawa, Poland}\email{nurowski@fuw.edu.pl}
%
\date{\today}
 \begin{abstract}
 We define pure radiation metrics with parallel rays  to be $n$-dimensional pseudo-Riemannian metrics that admit a parallel null line bundle $K$ and whose Ricci tensor vanishes on vectors that are orthogonal to $K$.
We give necessary conditions in terms of the Weyl, Cotton and Bach tensors for a pseudo-Riemannian metric to be conformal to a pure radiation metric with parallel rays. Then we derive 
  conditions in terms of tractor calculus that are equivalent to the existence of a pure radiation metric with parallel rays in a conformal class. We also give an analogous result for n-dimensional pseudo-Riemannian  pp-waves.
\\[.3cm]
\end{abstract}

\maketitle

\section{Introduction}
In general relativity, i.e. on a four dimensional Lorentzian manifold identified as space-time,  an energy momentum-tensor $T_{ab}$ satisfying the conditions
\begin{equation}\label{grpurerad}
T_{ab}=\phi^2 K_aK_b\ \text{ with a null vector }K_aK^a,
\end{equation}
is called an energy momentum tensor of a {\em pure radiation} or of a {\em null dust} (see \cite{stephani} for an overview). Ignoring the sign in \eqref{grpurerad}, this is equivalent to the existence of a null vector $K^a$ such that
\[
T_{ab}X^a=0
\text{ for all $X^a$  with } K_aX^a=0.\]
In case of a zero cosmological constant, 
this equation, 
via the Einstein field equations, 
\[
R_{ab}-\tfrac{1}{2}Rg_{ab}=T_{ab},
\]
becomes 
\begin{equation}\label{ricpurerad}
R_{ab}X^a=0
\text{ for all $X^a$ with } K_aX^a=0.\end{equation}
Here by
 $R_{ab}$ w denote the Ricci tensor of the metric. 

Another condition that appears in general relativity is the existence of a 
null vector $K^a$ that spans a parallel ray, i.e.
\begin{equation}
\nabla_aK^b=f_aK^b.
\end{equation}

Such metrics belong to the 
 so-called Kundt's class and have
 special Lorentzian holonomy. 
 If in addition $\nabla_{[a}f_{b]}\equiv 0$, the vector $K^a$ can be rescaled to a parallel null vector. 
 A particularly interesting class of metric in Kundt's class are those for which we have a {\em parallel null vector} and the energy momentum tensor is of pure radiation. Such metrics are called pp-waves in general relativity.

In the present article we are interested in generalisations of pure radiation metrics with parallel rays  to arbitrary dimension and signature and we will focus on their conformal properties. First, let us define the class of metrics we are interested in. 
\bde \label{defn1}
 Let 
 $g$ be a pseudo-Riemannian metric on an $n$-dimensional manifold.
\bnum
\item[(i)] \label{null-ric-def} The metric $g$ is a {\em  pure radiation metric}
if there is a null vector field $K^a$ such that
\begin{equation*} R_{ab}=\phi K_aK_b,\end{equation*}
with a function $\phi$ and $R_{ab}$ being the Ricci tensor of the metric.
\item[(ii)]\label{raydef} The metric $g$ has {\em parallel rays}
if there is a null vector field $K^a$ such that 
 $K^a$ spans  a  parallel null line bundle,  i.e.
\begin{equation*} \nabla_aK^b=f_aK^b
\end{equation*}
\item[(iii)] If $g$ satisfies both conditions (i) and (ii) for the same vector $K^a$ it is called {\em pure radiation metric with parallel rays} or {\em aligned pure radiation metric}.
\enum
\ede

The property of a metric to be a pure radiation metric is equivalent to 
\begin{equation}\label{null-ric}
R_{ab}X^a=0
\text{ for all $X^a$ with } K_aX^a=0.\end{equation}
Hence, pure radiation metrics have vanishing scalar curvature. This implies that the curvature $R_{abcd}$ and the Weyl tensor   $C_{abcd}$ of a pure radiation metrics satisfy
\begin{equation}
\label{no1}
R_{abcd}X^aY^b=C_{abcd}X^aY^b
\end{equation}
for all $X^a$ and $Y^b$ orthogonal to $K^a$.
Note also that that the Ricci tensor satisfies property \eqref{null-ric} if and only if the Schouten satisfies  property \eqref{null-ric}.

Note that for a metric with parallel rays spanned by the vector field $K^a$, i.e. with 
\begin{equation}\label{ray}
\nabla_aK^b=f_aK^b,
\end{equation}
the vector $K^a$ can be rescaled to a parallel vector field if and only if
\[\nabla_{[a}f_{b]}=0.\]
Furthermore, one can show  that for a metric with parallel rays defined as  in (ii) of Definition~\ref{defn1} we always find a null vector $K^a$ spanning the rays, such that
\[\nabla_a{K}_b= \psi{K}_a{K}_b,\]
with a function $\psi$.


In the first part of the paper we will derive necessary conditions in terms of the Weyl, Cotton and Bach  tensors for a pseudo-Riemannian metric to be {\em conformal}, i.e. locally conformally equivalent,  to a pure radiation metric with parallel rays. Then, as a special class of such metrics,  we study the  pp--waves in arbitrary signature. In the main part of the paper we use the normal conformal tractor calculus of \cite{bailey-eastwood-gover94} in order to derive {\em equivalent} conditions for a conformal class to contain a pure radiation metric with parallel rays. In Theorem \ref{theo-null} we prove: {\em A conformal class contains a pure radiation metric with parallel rays  if and only if the tractor bundle contains a sub-bundle of totally null 2-planes that is parallel with respect to the normal conformal tractor connection.} 
As a corollary we obtain a characterisation of conformal pp--waves in terms of the tractor connection.

\section{Tensorial obstructions}
In this section we derive tensorial obstructions for a metric to be conformal to an aligned pure radiation metric. Our conventions are as in \cite{gover/nurowski04} with the basic tensors derived from the pseudo-Riemannian metric $g_{ab}$ defined as follows: the curvature tensor, ${R_{ab\  d}^{\ \ c}}$ defined by
$(\nabla_a\nabla_b-\nabla_b\nabla_a)X^c=R_{ab\  d}^{\ \ c} X^d 
$ and its trace, the Ricci tensor, defined by $R_{ab}:=R_{cb\  d}^{\ \ c}$ and the scalar curvature, defined by$R:=R_a^{\ a}$. 
Fundamental tensors in conformal geometry are the Schouten tensor
\[\ro_{ab}:=\tfrac{1}{n-2}\left(R_{ab}-\tfrac{R}{2(n-1)} g_{ab}\right),\]
the Weyl tensor
\[C_{abcd}:=R_{abcd}-
2(g_{c[a}\ro_{b]d}+g_{d[b}\ro_{a]c}),
\]
as well as the Cotton and the Bach tensor
\be 
A_{abc}&:=&2\nabla_{[b}\ro_{c]a}\\
 B_{ab}&:=&\nabla^c
A_{acb}+\ro^{dc}C_{dacb}.
\ee
The Cotton tensor has a cyclic symmetry and satisfies
\be
 (n-3)A_{abc}&=&\nabla^d
C_{dabc}
\\
 \nabla^a A_{abc}&=&0.
\ee
Now we derive properties of these tensors when the metric is an aligned pure radiation metric.
\bs\label{pureprop}
Let $g$ be a pure radiation metric with parallel rays spanned by  $K^a$.
Then the Weyl tensor $C$, the Cotton tensor $A$ and the Bach tensor $B$ of $g$ satisfy
\begin{eqnarray}
C_{abcd}K^aX^b&=&0\label{weyl0}
\\
A_{abc}X^a&=&0\label{cotton0}
\\
B_{ab}X^a&=&0
\label{bach0}
\end{eqnarray}
for all $X^a$ orthogonal to $K^b$, i.e. with $X^aK_a=0$.
\es

\begin{proof}
Pure radiation implies  $R=0$ and we get  for the Schouten tensor $\ro$ that
\beq \label{purerad}
\ro_{ab}=\phi K_aK_b, 
\eeq
for a function $\phi$.
For $X^a$ and $Y^b$ orthogonal to  $K^a$  equation \eqref{purerad} gives
\[
C_{abcd}X^cY^d= R_{abcd}X^cY^d 
-2\phi \left(  g_{c[a}K_{b]}K_{ d} + g_{d[b}K_{a]}K_{c}\right)X^cY^d
= R_{abcd}X^cY^d.
\]
Equation \eqref{ray} gives
\[
R_{abcd}K^d= \nabla_{[a} f_{b]}K_c\]
which implies
$ R_{abcd}X^c K^d=0$
for all $X^a$ orthogonal to $K^a$.
Hence, $ C_{abcd}X^c K^d=0$.

Equation \eqref{cotton0} follows immediately from the definition of $A_{abc}$, from $R_{ab}X^b=\ro_{ab}X^b=0$ and the fact that $K^\bot$ is parallel, i.e. $\nabla_a X^bK_b=0$. 

Finally, we prove \eqref{bach0}. By the definition of the Bach tensor and from equations \eqref{cotton0}, \eqref{weyl0} and  \eqref{purerad} we get
\[
B_{ab}X^a= \nabla^cA_{acb}X^a+\ro^{cd}C_{acbd}X^a
=
\phi C_{acbd}X^aK^cK^d=0
\]
for all $X^b$ orthogonal to $K^b$. By the symmetry of $B$, this proves \eqref{bach0}. 
\end{proof}
Now we will use this proposition in order to derive obstructions for a metric to be conformal to an aligned  pure radiation metric. Under a conformal change of the metric,
$\hat g_{ab}=\e^{2\ups}g_{ab}$, the Levi-Civita changes as
\beq
\label{lc}
\hat\nabla_aX_b\ =\ \nabla_aX_b-\ups_aX_b -\ups_bX_a+g_{ab}\ups^dX_d,\eeq
with $\ups_a=\nabla_a\ups$. The Schouten tensor transforms as
\beq\label{schouten}\hat\ro_{ab}\ = \ \ro_{ab}-\nabla_a\ups_b+\ups_a\ups_b -\einhalb \ups_c\ups^cg_{ab}.\eeq
Using this, we can show:
\btheo\label{theo1}
Let $g_{ab}$ be a metric that is conformal to a  pure radiation metric with parallel rays $K^a$. Then the Weyl tensor of $g_{ab}$ satisfies
\beq\label{weyl}
C_{abcd}K^cX^d\ =\ 0\, ,\ \text{ for all $X^d$ orthogonal to $K^a$, i.e.  with $K_dX^d=0$.	}
\eeq 
Furthermore, 
there is a gradient field $\ups^a$  such that
\begin{eqnarray}
\label{cotton}
\left( A_{abc} + C_{cba}^{\ \ \ d}\ups_d\right)X^a&=& 0 \\
\label{bach}
\left(B_{ab}- (n-4) C_{acbd}\ups^c\ups^d\right) X^aY^b&=&0
\end{eqnarray}
for all $X^a$ and $Y^b$ orthogonal to $K^a$.
\etheo
\bbem In four dimensions and Lorentzian signature, condition \eqref{weyl} means that $W$ is of Petrov type III or N.
\ebem
\bprf
Let $\hat g=\e^{2\ups}g$ be a pure radiation metric with parallel  rays spanned by $K^a$ and with $\hat R_{ab}X^b=0$ for all $X^a$ orthogonal to $K^a$.
Then equations \eqref{weyl} and \eqref{cotton} follow immediately from the proposition by the conformal invariance of the Weyl tensor and the conformal transformation formula for the Cotton tensor:
\[ \hat A_{abc}= A_{abc}+ C_{cba}^{\ \ \ d}\ups_d.\]
In order to prove equation \eqref{bach}, we take the divergence of  \eqref{cotton}.  Using the definition of the Bach tensor,
the well known identity for the divergence of the Weyl tensor
\begin{equation}\label{divweyl}
\nabla^b C_{cbad}=(n-3)A_{cad},
\end{equation}
 and the transformation formula \eqref{schouten} in order to eliminate the $\nabla^c\ups_d$-terms, we get
 \begin{eqnarray}\nonumber
 0&=&\left(\nabla^bA_{abc}+\nabla^b(C_{cba}^{\ \ \ d}\ups_d)\right)X^a +
 \left( A_{abc} + C_{cba}^{\ \ \ d}\ups_d\right)\nabla^bX^a
 \\
 &=&
 \left( B_{ac} -(n-3) A_{acd}\ups^d + C_{cbad}\left( \ups^b\ups^d -\hat\ro^{bd}\right)\right)X^a
 \\
 &&{ }+
  \left( A_{abc} + C_{cba}^{\ \ \ d}\ups_d\right)\nabla^bX^a
  \nonumber
\label{add}
  \end{eqnarray}
 By equations \eqref{purerad} and \eqref{weyl} we get 
 \[ 
  C_{acbd}\hat\ro^{cd}X^b = \phi \,  C_{acbd}K^cK^dX^b =0.\]
Since only  $\hat\nabla^bX^a$ is orthogonal to $K^a$ but not necessarily $\nabla^bX^a$ we compute  using \eqref{lc} 
\be
\left( A_{abc} + C_{cba}^{\ \ \ d}\ups_d\right)\nabla^bX^a
&=&
\left( A_{abc} + C_{cba}^{\ \ \ d}\ups_d\right) \ups^aX^b
\\
&=&
-(A_{bca}+A_{cab})\ups^aX^b
\\
&=&
(C_{acb}^{\ \ \ d}\ups_d-A_{cab})\ups^aX^b
\\
&=&
-(C_{bac}^{\ \ \ d}\ups_d+A_{cab})\ups^aX^b
\ee
since both, the Cotton and the Weyl tensor are trace free and satisfy the Bianchi identity.
Hence, when contracting \eeqref{add} with $Y^c$ we get on one hand
\[
\left( A_{abc} + C_{cba}^{\ \ \ d}\ups_d\right)(\nabla^bX^a)Y^c=0,\]
and on the other hand,  by  \eqref{cotton},
 \[
 A_{acd}\ups^dX^aY^c
 =
 -C_{dcab}\ups^d\ups^bX^aY^c
 =
 C_{cdab}\ups^b\ups^dX^aY^c
 \]
 for $X^b$ and $Y^a$ orthogonal to $K^a$. This yields 
 equation \eqref{bach}.
\eprf
Finally, we return to metrics with parallel rays. From the transformation formula \eqref{lc} it is obvious that the property of admitting  parallel rays is not invariant under a general conformal change of the metric.
The following proposition shows that a parallel rays cannot be improved to a parallel vector field under a conformal change.
\bs
Let $g$ be a pseudo-Riemannian metric with parallel rays spanned by  $K^a$ which cannot be rescaled to a parallel vector. Then there is also no metric in the conformal class of $g$ for which 
$K^a$ can be rescaled to a parallel vector field.\es
\bprf This statement follows from \eqref{lc}: Assume that $K^a$ spans a parallel ray, i.e. $\nabla_aK_a=f_aK_a$ for $K_a=g_{ab}K^b$, such that it cannot be rescaled to a parallel vector field. This means that $\nabla_{[a}f_{b]}\not=0$. For a conformally changed metric 
$\hat g_{ab}=\e^{2\ups}g_{ab}$
the metric dual of $K^a$ with respect to $\hat g$ is $\hat K_a =\e^{2\ups}K_a$. Then, that  the rescaled one-form $\vf\hat K_a$, where $\vf$ is a function $\vf$ without zeros, is parallel, reads as

\be
0&=& 
\e^{2\ups}\left( \nabla_a\vf K_b+ 2 \vf \ups_a +\vf \hat\nabla_aK_b\right)
\\
&=&
\e^{2\ups}\left( \nabla_a\vf + \vf \ups_a+\vf f_a\right) K_b+ \e^{2\ups}\vf \left( g_{ab}\ups^dK_d- \ups_bK_a\right).
 \ee
Contracting this with arbitrary $X^a$ and $Y^b$ from $K^\bot$ shows that $\ups^dK_d=0$. Then contracting only with $X^b$ gives $\ups_bX^b=0$ which shows that $\ups_b=\psi K_b$, and therefore $\ups_bK_a=\ups_aK_b$  Hence, the above equation gives
\[
0= 
\nabla_a\vf +\vf f_a=\vf (\nabla_a \log \vf+f_a)
\]
which contradicts $\nabla_{[a}f_{b]}\not=0$.
\eprf

\section{pp--waves as aligned pure radiation metrics}
Now we will derive  similar conditions for a more special class of aligned pure radiation metrics that can be considered as pp-waves in arbitrary signature. For a given null vector $K^a$, not necessarily parallel, we will consider the  
fundamental property of the curvature tensor $R$, 
\begin{equation}
\label{pp}
R_{abcd}X^aY^b=0\ \text{ for all $X^a$ and $Y^b$ orthogonal to $K^a$.}
\end{equation}
Then, in analogy with the Lorentzian case, we define:
\bde
An $n$-dimensional  pseudo-Riemannian metric is a {\em pp--wave} if it admits a parallel null vector field 
for which property \eqref{pp} is satisfied.
\ede
These metrics have a highly degenerated curvature tensor and hence a very special restricted holonomy group. By the restricted holonomy group we refer to the group of parallel transports along contractible loops. The restricted holonomy group  is the connected component of the full holonomy group.  
In analogy to the Lorentzian case, if the signature of the metric is $(p,q)$, one can show that these metrics have their restricted holonomy groups contained in the abelian normal subgroup $\rr^{p-1,q-1}$ in the group $\SO^0(p-1,q-1)\ltimes \rr^{p-1,q-1}$, where the latter is the  is the stabiliser in $\SO^0 (p,q)$ of a null vector. In particular, if we denote by $\cal K$ the line distribution given by $K^a$, then the vector bundle $S:=\cal K^\bot /\cal K$ is flat with respect to the connection $\nabla^S_X[Y]:=[\nabla_XY]$ for $Y$ a section of $\cal K^\bot$.

Again in analogy to the Lorentzian case \cite{leistner-nurowski08},  the conformal Fefferman-Graham ambient metric and the tensor obstructing its existence can be computed explicitly for pp-waves.

In the following we will describe how pp-waves relate to pure radiation metrics.
First, we verify that pp--waves  are  pure radiation metrics with parallel rays spanned by the parallel null vector field:
\bs\label{ppprop}
Let $g$ be a pseudo-Riemanian pp--wave metric with a parallel null vector $K^a$.
Then $g$ is a pure radiation metric, i.e. equations \eqref{cotton0} and \eqref{bach0} are satisfied, and furthermore we have
\begin{eqnarray}
C_{abcd}K^a&=&0
\label{ppN}\\
\label{ppweyl}
C_{abcd}X^aY^b&=&0\\
 A_{cab}X^aY^b&=&0,
 \label{ppcotton}
 \end{eqnarray}
 for all $X^a$ and $Y^b$ orthogonal to $K^a$.
\es
\bprf
Since $K^a$ is parallel and because of property \eqref{pp} we get $R_{ab}X^a=0$ for all $X^a$ orthogonal to $K^a$. Hence, $R_a^{\ a}=0$ which implies $\ro_{ab}X^a=0$ and $R_{abcd}X^aY^b=C_{abcd}X^aY^b$ for all $X^a$ and $Y^a$ orthogonal to $K^a$.
Then the required \eeqref{ppcotton} immediately follows from the identity \eqref{divweyl} for the divergence of the Weyl tensor.
Finally, as $K^a$ is parallel, it follows
\be
0&=& R_{abcd}K^a
\\
& = &C_{abcd}K^a
+2(g_{c[a}\ro_{b]d}+g_{d[b}\ro_{a]c})K^a
\\
&=&
C_{abcd}K^a
+2K_c\ro_{bd} -2g_{cb}\ro_{ad}K^a
+2g_{db}\ro_{ac}K^a- 2K_d\ro_{bc}
\\
&=&
C_{abcd}K^a
+2K_c\ro_{bd} 
- 2K_d\ro_{bc}.
\ee
which shows that $C_{abcd}K^aX^c=0$ for all $X^c$ orthogonal to $K^a$. But this implies that 
$C_{abcd}K^a=0$.
\eprf
The following proposition shows that for pure radiation metrics with parallel rays the property \eqref{pp} implies that the rays contains a parallel null vector. Recall \eqref{no1} which implies that a pure radiation metric satisfies 
 \eqref{pp} for the curvature tensor if and only if it satisfies  \eqref{ppweyl} for the Weyl tensor.
\bs\label{ppweylprop}
A pure radiation metric with parallel rays $K^a$ that satisfies property  \eqref{pp}, or the equivalent property \eqref{ppweyl},  admits a parallel null vector field in direction of $K^a$ and hence, is a pp--wave.
\es
\bprf
Assume that $\nabla_aK^b=f_aK^b$ with a one-form $f_a$. Then, locally $K^a$ can be rescaled to a parallel vector field if and only if $f_a$ is closed, i.e. iff $\nabla_{[a}f_{b]}=0$. Since
\[R_{abc}^{\ \ \ d}K^c= \nabla_{[a}f_{b]}K^d,\]
we have to show that $R_{abc}^{\ \ \ d}K^c=0$. Now \eqref{pp} implies that
\[
0\ =\  R_{abc}^{\ \ \ d}X^aY^bK^c\ =\ R_{abc}^{\ \ \ d}X^aK^b
\] for all $X^a$ and $Y^b$ orthogonal to $K^a$. Hence, the only possibly non-vanishing terms of $R_{abc}^{\ \ \ d}K^c$ could be
$R_{abcd}L^aX^bK^cL^d$ with $L^a$ a vector transversal to $K^\bot$ and $X^b\in K^\bot$. But since we have pure radiation, also these terms vanish
\[0= R_{ab}X^aL^b=R_{cabd}L^cX^aL^b K^d.\]
Indeed, if 
$(K,L,E_1, \ldots , E_{n-2})$ is a basis with $K^\bot$ spanned by $(K, E_1, \ldots , E_{n-2})$, $L^\bot$ is spanned by $(L, E_1, \ldots , E_{n-2})$, and such that $g(K,L)=1$ and $g(E_i,E_j)=\pm\delta_{ij}$, we have
\[0=Ric(Y,L)=R(L,Y,K,L)+\sum_{i=1}^{n-2}\eps_i R(E_i,Y,E_i,L)= R(L,Y,K,L).\]
Here $Ric(.,.)$ denotes the Ricci tensor and $R(.,.,.,.)$ the $(4,0)$-curvature tensor.
\eprf
%

We also have the following proposition.
\bs
Let $(M,g)$ be pseudo-Riemannian manifold of dimension $n>2$ with parallel null vector field $K^a$, i.e. with $\nabla_aK^b=0$, satisfying \eqref{ppweyl} for the Weyl tensor. Then $g$ is a pure radiation metric aligned with  $K^a$ and hence a pp--wave.
\es
\bprf
First we show that $g$ has vanishing scalar curvature $R=0$.
Since $K^a$ is parallel we have that $R_{ab}K^a=0$ and $R_{abcd}K^a=0$. Because of property \eqref{ppweyl}, with the same notation and same  basis $(K,L, E_1, \ldots , E_{n-2})$ as in the proof of previous proposition we get for the scalar curvature 
\be
R&=&
\sumij \eps_i\eps_jR(E_i, E_j, E_i, E_j)
\\
&=&
2(n-3)\sumi \ro\eps_i (E_i,E_i)
\\
&=&
2(n-3)\left( \trace (\ro ) -\ro (K,L)\right)
\\
&=&
2(n-3)\left( \frac{R}{2(n-1)}  + \frac{R}{(n-1)(n-2)}g(K,L)\right)
\\
&=&
\frac{n(n-3)}{(n-1)(n-2)} R,
\ee
which shows that $R=0$.
Next, we verify that the metric is a pure radiation metric. For $X\in K^\bot $ and $V\in TM$, since $K$ is parallel and the Weyl tensor satisfies \eqref{ppweyl}, we get
\be
Ric(X,V)&=& \sumi\ \eps_i R(E_i,X,E_i,V)
\\
&=&(n-2)\ro(X,V)- \sumi \eps_i(g(E_i,X)\ro(E_i,V)+g(V,E_i)\ro(X,E_i)),
\ee
in which we have used for the last equality that $R=0$. Hence
we get
\[0=\sumi \eps_i(g(E_i,X)\ro(E_i,V)+g(V,E_i)\ro(X,E_i))
=
2\ro(X,V)-g(K,V)\ro(L,X),
\]
which shows that $\ro(X,V)=0$. Hence, $g$ is a pure radiation metric and thus, by \eqref{no1}, a pp-wave.
\eprf 

 From Theorem \ref{theo1} we get the following obstruction for the existence of a pp--wave metric in a conformal class:

\bfolg
Let $g$ be a metric that is conformal to a  pp-wave metric with parallel null vector $K^a$.
Then, in addition to properties \eqref{cotton} and \eqref{bach}, for  the Weyl and Cotton tensors   of $g$ 
we have equations \eqref{ppN}, \eqref{ppweyl} and \eqref{ppcotton}.
\efolg
\bprf
Let $\hat g=\e^{2\ups}g$ be  a pp--wave metric with parallel null vector
 $K^a$.  
 First note that the rays spanned by $K^a$ and its orthogonal complement are conformally invariantly defined. Both are parallel with respect to $\hat\nabla$ but not parallel for $\nabla$.
 
By Proposition \ref{ppprop}, we have that $\hat{C}_{abcd}K^a=0$ and  $\hat C_{abcd}X^aY^b=0$  for all $X^a$ and $Y^b$ orthogonal to $K^a$. The conformal invariance of the Weyl tensor then implies equations \eqref{ppN} and \eqref{ppweyl}. 

In order to derive equation \eqref{ppcotton} we apply Theorem \ref{theo1}.
The Bianchi identity for the Cotton and Weyl tensors together with \eqref{cotton}  for $X^a$ and $Y^b$ orthogonal to $K^a$
gives
\be
0&=&
(A_{abc}+A_{bca}+A_{cab})X^bY^c
\\
&=&
A_{abc}X^bX^c
-  C_{acb}^{\ \ \ d}\ups_d X^bY^c
-  C_{bac}^{\ \ \ d}\ups_d X^bY^c
\\
&=&
A_{abc}X^bX^c
+  (C_{cba}^{\ \ \ d}+ C_{bac}^{\ \ \ d})\ups_d X^bY^c
-  C_{bac}^{\ \ \ d}\ups_d X^bY^c.
\ee
Then \eqref{ppweyl} is $ C_{bca}^{\ \ \ d} X^bY^c=0$ and implies the required \eeqref{ppcotton}.
\eprf

\section{Conformal standard tractors}

An invariant way of describing the geometry of a conformal manifold $(M,[g])$ is provided by the {\em normal conformal Cartan connection} $\w$ (refer to \cite{kobayashi72}, \cite{bailey-eastwood-gover94}, \cite{cap/gover02} or to the survey \cite{baum07} for details of the following). For example, $\w$ detects local Einstein metrics in the conformal class, which correspond to covariantly constant sections with respect to $\w$. Now we will establish a condition on $\w$ that is equivalent to the existence of a local pure radiation metric  in the conformal class. 

The normal conformal Cartan conection 
 is defined uniquely by the conformal class $[g]$ and defined by the following data:
Let $(p,q)$ be the signature of the conformal class and $P\subset \SO_0(p+1,q+1)$ the stabiliser in $\SO_0(p+1,q+1)$  of a null line $I$. Then $\w$ is a Cartan connection with values in $\laso(p+1,q+1)$ on a principal  $P$-bundle $\cal P$, the {\em Cartan bundle}. Furthermore, it satisfies  a normalisation condition that makes it unique.  
 As a Cartan connection, $\w$ defines an invariant absolute parallelism and hence, gives no horizontal distribution in $T\cal P$. 
In order to use the usual principle fibre bundle formalism, one extends the  Cartan bundle to  a bundle $\ol{\cal P}$ with structure group $\mathrm{SO}^0(p+1,q+1)$ via $\ol{\cal P}=
\cal P\times_P \mathrm{SO}^0(p+1,q+1)$ on which the Cartan connection $\w$ extends to a principle fibre bundle connection $\ol{\w}$.

By associating the standard representation $\rr^{p+1,q+1}$ of $\mathrm{SO}(p+1,q+1)$ to $\ol{\cal{P}}$ we obtain a vector bundle $\cal T$ of rank $p+q+2$, called 
{\em standard tractor bundle}. 
${\cal T}$ is equipped
with the covariant derivative $\ol{\nabla}$ induced by 
 $\ol{w}$. Since $\ol{\w}$ is an $\laso(p+1,q+1)$-connection,
  there is an invariant metric $h$ on $\cal T$. The triple $(\cal T, \ol{\nabla}, h)$ 
 is called {\em normal conformal standard tractor bundle}. 
 
 Let $I\subset \rr^{p+1,q+1}$ be the null line defining $P$, i.e. $P$ is the stabiliser  group of $I$. Then, corresponding to the filtration
 $I\subset I^\bot \subset \rr^{p+1,q+1}$, there is a bundle filtration of the tractor bundle
 \[ \cal I\subset \cal I^\bot\subset \cal T\]
 into a null line bundle and its orthogonal complement, w.r.t. to the tractor metric. Then the projection $p:\cal P\to \cal P^0$  of the Cartan bundle onto the bundle of conformal frames $\cal P^0$ defines the following projection 
\[ \begin{array}{rcccl}
 \mathrm{pr}_{TM}\ :\  \cal I^\bot=\cal P\times_P I^\bot & \to & \cal I^\bot/\cal I & \simeq & TM \ \simeq \ \cal P^0\times_{\mathrm{CO}_0(p,q)} \left( I^\bot/I\right)  \\
 \left[\vf, v\right] &\mapsto & \left[\vf , [v] \right]&\mapsto & \left[p(\vf),[v]\right]
 \end{array}
 \]
 in an invariant way.
Via this projection, the tractor metric defines the conformal structure.
 
Every metric $g$ in the  
 conformal class $[g]$ induces a splitting 
\[\left(\cal T,   h\right)\ \simeq \underline{\rr}\+ TM\+ \underline{\rr},\]
where each $\underline{\rr}$ denotes the trivial line bundle $M\times \rr$. 
Changing the metric in the conformal class, $\hat g_{ab}=\e^{2\ups}g_{ab}$, this splitting transforms via
\begin{equation}\label{splitting}
\left( \begin{array}{c} 
\hat \rho  \\ \hat X^a \\ \hat\s \end{array}\right)
=
\left(
\begin{array}{c}
\rho - \ups_aX^a - \einhalb \s \ups^b\ups_b 
\\ X^a+ \s\  \ups_a  \\
\s 
  \end{array}\right).
\end{equation}
In every such splitting the line bundle $\cal I$ is spanned by
\[I^A:= \left(
\begin{array}{c} 1\\0\\0\end{array}\right),\]
where capital indices range in $0,1,\ldots ,n , n+1$, with $n$ the dimension of $M$.
 The  bundle metric $h$ of signature $(p+1,q+1)$ on $\cal T$ is given by
\[
h\left( \left( \begin{array}{c} \rho \\ X^a \\ \s \end{array}\right), \left( \begin{array}{c} \alpha \\Y^b \\\beta\end{array}\right)\right)
=
\s\alpha+\rho\beta +g_{ab}X^aY^b. \]
Hence, 
 $\cal I^\bot$ is given by tractors of the form
 \[\left( \begin{array}{c} \rho \\ X^a \\ 0 \end{array}\right)\]
 with tangent vectors $X^a$,
and that the projection $\mathrm{pr}_{TM}$ is invariantly given as
\[\mathrm{pr}_{TM} : \cal I^\bot \ni \left( \begin{array}{c} \rho \\ X^a \\ 0 \end{array}\right) \mapsto X^a\in TM.\]

In a splitting defined by a metric $g_{ab}$ the tractor connection  $\ol{\nabla}$ is given by the following formula,
\begin{eqnarray}\label{nabtrac}
\ol{\nabla}_a\left( \begin{array}{c} \rho \\X^b \\[-1mm] \s \end{array}\right)
&=& 
\left(\begin{array}{c} 
\nabla_a \rho-\ro_{ab}X^b
\\
 \nabla_aX^b +\rho \delta_a^b +\s \ro_{a}^{\ b} 
\\
\nabla_a\s -g_{ab}X^b
\end{array}\right),
\end{eqnarray}
where $\nabla $ is the Levi-Civita connection and $\ro$ the Schouten tensor of the metric $g$. The curvature of this connection is given as
\begin{equation}
\label{tracurv}
\ol{R}_{ab}
\left( \begin{array}{c} \rho \\X^c \\[-1mm] \s \end{array}\right)
=
\left( \begin{array}{c} 
-A_{dab}X^d
 \\ \s A^c_{\ ab}+C_{ab\  d}^{\ \ c}X^d \\ 0 \end{array}\right).
\end{equation}
It is a well known fact, that covariantly constant sections of $\ol{\nabla}$ are in one-to-one correspondence with local Einstein metrics on an open and dense subset of $M$. In fact, if 
$X^B=
\left( \begin{array}{c} \rho \\X^b\\[-1mm]\s \end{array}\right)$ is covariantly constant for $\ol{\nabla}$ if and only if  $\s^{-2} g$ is an Einstein metric on the open and dense complement of the zero set of $\s=X^AI_A$. 
Indeed, $X^B$ is parallel if and only if $X^b=\nabla^b\s$ and $\nabla_a\nabla_b \s +\s P_{ab}+\rho g_{ab}=0$, where the latter is equivalent to $\s^{-2}g$ being an Einstein metric outside the zero set of $\s$. 
Furthermore, applying the tractor curvature $\ol{R}_{ab}$ to a parallel tractor $X^B$ we get
the integrability conditions $\ol{R}_{ab}X^C=0$, i.e. for $X_c=\nabla_c\s$ we get
\be
A^c_{\ ab}X_c&=&0
\\
\s A_{cab}+C_{abc}^{\ \ \  d}X_d&=&0.
\ee
Setting $\s=\e^{-\ups}$ we get $X_a=-\s \ups_a$ with $\ups_a=\nabla_a\ups$ and obtain 
\[
 A_{cab}-C_{abc}^{\ \ \  d}\ups_d=0,
 \] as one of the integrability conditions that was derived in \cite{gover/nurowski04}. 

\section{Tractorial characterisation of aligned pure radiation metrics}

Now we will describe the situation in which the normal conformal tractor holonomy admits an invariant plane. The case when the plane is non-degenerate implies that the conformal class contains a product of Einstein metrics with related Einstein constants (see \cite{armstrong07conf} for Riemannian conformal classes and related results  in the unpublished parts of  \cite{leitner04} for arbitrary signature). Here we will deal with the case that the null plane is totally null in arbitrary signature. The following theorem is a generalisation to arbitrary signature of the corresponding  result for Lorentzian conformal classes which we proved in \cite{leistner05a}.

\btheo\label{theo-null}
Let $(M,[g])$ be a pseudo-Riemannian conformal manifold of dimension $n>2$. Then the normal conformal Cartan connection admits a parallel totally null plane if and only if, on an open and dense subset of $M$,  there is an aligned pure radiation metric $g_{ab}$ in the conformal class  $[g]$.
\etheo

\bprf
Let $\cal T$ be the normal conformal tractor bundle.
First, we prove a Lemma.
\blem
 Given a bundle of null lines  in $TM$ spanned by $K^a$  and a metric $g_{ab}$ in the conformal class, then the following statements are equivalent:
\bnum
\item[(i)] $g_{ab}$ is pure radiation metric with parallel rays spanned by $K^a$,
\item[(ii)]
W.r.t. the splitting of $\cal T=\underline{\rr}\+TM\+\underline{\rr}$ given by $g_{ab}$, the plane bundle $\cal H$ spanned by
 \begin{equation}
\label{null-plane}
K^A=\left(\begin{array}{c}0\\ K^a\\0\end{array}\right)\text{ and }J^A=\left(\begin{array}{c}0\\0\\1\end{array}\right)
\end{equation}
is parallel for the tractor connection.
\enum
\elem
\bprf
Since the tractor connection preserves the metric, $\cal H$ is parallel if and only if \[\cal H^\bot =
\left\{\left(\begin{array}{c}0\\X^b\\
\s\end{array}\right)\mid K_aX^a=0, \s\in \rr\right\}\] is parallel. 
By pairing with $K^A$ and $J^A$ we get that $\cal H^\bot $ being parallel is equivalent to
\[
P_{ab}X^b=0\text{ and } \left(\nabla_aX^b-\s P_{a}^{\ b} \right)K_b=0,\]
which is equivalent to $R_{ab}X^b=0$ and $\nabla_aK^b=f_aK^b$.
\eprf
The 
 `if'-statement of the theorem follows immediately from the lemma:
If $g_{ab}$ is  an aligned pure radiation metric in the conformal class, we split the tractor bundle with respect to $g_{ab}$  and get a parallel null plane $\cal H$  for the tractor connection as in \eqref{null-plane}.

Now we show the `only if'-statement in the theorem. Assume that $\cal H$ is a   null plane bundle that is parallel for the tractor connection. We have to find a null vector $K^a$ and a metric in the conformal class such that $\cal H$ is spanned by $K^A$ and $J^A$ as in 
as in \eqref{null-plane}. 
First, using that $\cal H$ is parallel and some basic linear algebra, we will define  the null line. 
Let $\cal I\subset \cal T$ be the null line bundle defined by the conformal structure. Then we define
\[\cal L:=\cal I^\bot\cap \cal H =\{X^A\in \cal H\mid X^AI_A=0\}\]
a subbundle of $\cal H$. We have:
\blem
$\cal L$ is a null line bundle and, over an open and dense subset of $M$, its invariant projection $\mathrm{pr}_{TM}(\cal L) $ is not zero.
\elem
\bprf As $\cal H$ is bundle of 2-planes,  for dimensional reasons we have $\cal L\not=\{0\}$. Hence, the rank of $\cal L$ is either one or two.
Then
we fix {\em some} metric in the conformal class and w.r.t to the induced splitting 
$\cal T=\underline{\rr}\+TM\+\underline{\rr}$ of the tractor bundle  we get 
\[\cal L
= 
\left\{ L^A=\left(\begin{array}{c}\rho\\L^a \\0\end{array}\right)=\rho I^A+L^a \in \cal H\right\}.\]
If $\cal L$
had rank two over an open set $U$, we would get $\cal H|_U\subset (\underline{\rr}\+TM)|_U=\cal I^\bot|_U$. Hence, with $\cal H$ being parallel, the tractor derivative of every $L^A\in \cal H|_U$ is again in $\cal H|_U$ and hence in $\cal I$. Differentiating in any direction using formula \eqref{nabtrac}, over $U$ we get
\[0= \ol{\nabla}_b L^AI_A= -g_{ba}L^a.\]
But this means that $L^a=0$ which excludes then rank of $\cal L|_U$ being two. Hence there is no open set over which $\cal L$ has rank two. Therefore it must have rank one over all of $M$ and hence,
%
 $\cal L$ is a line bundle.

Furthermore, assume that 
$\mathrm{pr}_{TM}\cal L\{0\}$ over an open set $U$ of $M$. Then, over $U$, 
every section of $\cal L$  would be of the form
 $\rho I^A$ with $\rho \in C^\infty(U)$.  Differentiating yields
 \[
\ol{\nabla}_b(\rho I^A)=  \nabla_b\rho\, I^A + \rho \nabla_bI^A = \nabla_b\rho I^A + \rho \delta_b^{\ a}.
\]
Since $ \cal L\subset \cal H$ we get that $\ol{\nabla}_b(\rho I^A)= \nabla_b\rho I^A + \rho \delta_b^{\ a}$
is again in $\cal H$
because $\cal H$ was parallel. Since $U$ was open, we can differentiate in any direction of $TM$, which yields acontradiction to $\cal H$ being a plane.
\eprf
The lemma shows that, over an open and dense subset of $M$,  we can define a  line bundle in $TM$ as
\[L:=\mathrm{pr}_{TM}\cal L\]
in an invariant way. In the following we will restrict our computations to this open and dense subset without explicitly mentioning it again. Since $\cal H$ and hence $\cal L$ were totally null, $L$ is a bundle of null lines. Also, its orthogonal complement, $L^{\bot}\subset TM$, is given as a projection from $\cal T$:
\blem\label{orthlem}
The orthogonal complement of $L$ satisfies $ L^\bot=\mathrm{pr}_{TM} \left(\cal L^\bot\cap\cal I^\bot\right)$.
\elem
\bprf Let $\cal L$ be spanned locally by $I^{A}+K^{a}$ with a null vector field $K^a$.  If $\rho I^{A}+ X^{a}\in \cal L^{\bot}\cap\cal  I^{\bot}$, then $X^{b}K_{b}=0$, i.e. $\mathrm{pr}_{TM} \left(\cal L^\bot\cap\cal I^\bot\right)\subset L^{\bot}$. On the other hand, if $X^{a}\in L^{\bot}$, then $I^{A}+ X^{a}\in\cal L^{\bot}\cap \cal I^{\bot}$.
\eprf

Next, we prove 
\blem\label{horth}
$\cal L^\bot=\cal H^\bot\+ \cal I$
\elem
\bprf First note that
$\cal L^\bot$ contains $I^A$. On the other hand,  $\cal L\subset \cal H$ gives  $\cal H^\bot\subset \cal L^\bot$. 
Recalling that $\cal L = \cal H \cap \cal I^\bot$ is a line in the plane
 $\cal H$ shows that there is an element  $X^A\in \cal H$ such that $X^AI_A=1$.
 This implies 
 $\cal H^\bot\cap \cal I=\{0\}$ and counting dimensions completes the proof.
 \eprf

Then Lemmas \ref{orthlem} and 
 \ref{horth}
 imply that
\begin{equation}\label{lorth}
L^\bot= \mathrm{pr}_{TM}\left(\cal H^\bot\cap \cal I^\bot\right),\end{equation}
which provides us with 
\blem
The hyperplane bundle $L^\bot\subset TM$ is integrable.
\elem
\bprf
Let $X^{a}$ and $Y^{b}$ two local sections of $L^\bot$ such that $X_{a}Y^{a}=0$. Then there is a smooth function $\rho$ such that
$\rho I^{a}+ Y^{b}\in \cal H^\bot\cap \cal I^{\bot} \subset \cal H^\bot$. Since $\cal H^\bot $ is parallel, we get that
\[
X^a \ol{\nabla}_a \left( \begin{array}{c} \rho \\ Y^b \\ 0 \end{array}\right)
= \left(\begin{array}{c} X^a(\nabla_a \rho-\ro_{ab}Y^b)\\ X^a\nabla_a Y^b+\rho X^b \\ 0\end{array}\right)\in 
\Gamma(\cal H^\bot\cap \cal I^\bot),\]
By \eqref{lorth} this shows 
that $X^a\nabla_a Y^B \in \Gamma(L^\bot)$ for $X^a, Y^b\in \Gamma(L^\bot)$ orthogonal to each other. By fixing a basis of $\cal L^\bot$ that consists of mutually orthogonal vectors, the vanishing of the torsion of $\nabla$ implies that $L^\bot$ is integrable.
\eprf
Hence, so far, given the parallel plane distribution $\cal H$ in $\cal T$, we have invariantly constructed a null line bundle $L$ in $TM$ such that the $L^\bot$ is integrable.

The integrability of $ L^\bot$ allows us to define a second fundamental form for $L^\bot$ in the following way. From now on we fix a null vector field $K^a$ spanning $L$ and
define a bilinear form on $L^\bot$ by
\[\Pi_{ab} :=\nabla_a K^cg_{cb} |_{L^\bot\times L^\bot}.\]
Since $L^\bot $ is integrable, $\Pi^K$ is symmetric.
We define the trace of $\Pi^K$ as
\[H:=g^{ab}  \Pi_{ij}E^a_iE^b_i\in C^\infty (M)\]
where $E^a_1,\ldots , E^a_{n-2}$ are linearly independent in $L^\bot$ and $i,j$ range over $1, \ldots , n-2$. 
Since $K^a \Pi_{ab}=0$, this is independent of the chosen $E^a_i$'s.
Now we claim that there is a metric $\hat{g}_{ab}=\e^{2\Upsilon}g_{ab}$ in the conformal class such that the corresponding function $\hat{H}$ is zero. To this end we notice that the transformation formula for $\hat{\Pi}$ is given by
\[X^aY^b \hat{\Pi}_{ab}= X^aY^b \nabla_a \hat{g}_{cb}K^c
=
\e^{2\Upsilon}\left(\Pi_{ab} +K^c\ups_c g_{ab}\right)X^aY^b,
\]
for $X^a,Y^b\in L^\bot$. Hence,
\[\hat{H}=\e^{2\ups}\left( H+ (n-2) K^a\ups_a \right).\]
Now the differential equation 
\[K^a\ups_a=\frac{H}{n-2}\]
has always a solution $\ups$, which ensures that we can chose $\hat{g}_{ab}$ such that $\hat{H}\equiv 0$. Finally, to conclude the proof, we fix this metric, omit the hat, and split the tractor bundle with respect to the new metric. Now let  $\rho I^A+K^a$ be an arbitray section of $\cal L$. Since $\cal H$ is parallel, differentiating in direction $Y^b\in L^\bot$ yields that
\[
Y^b \ol{\nabla}_b\left( \begin{array}{c} \rho \\ K^a \\ 0 \end{array}\right)
= \left(\begin{array}{c} Y^b\nabla_b\rho-\ro_{ab}Y^bK^a
\\ Y^b \nabla_b K^a+\rho Y^a \\ 0\end{array}\right)\]
is still a section of $\cal H$, but also of 
 $\cal L$ since the last component vanishes. 
 This means that $Y^b\nabla_b K^a+\rho Y^a $ is in $L$ which implies that
\[\Pi_{ab} Y^aY^b+\rho g_{ab}Y^aY^b=0\]
for all $Y^a\in L^\bot$. Taking the trace, $H=0$ gives
$0=\rho (n-2)$,
which results in $\cal L=\rr\cdot\left(\begin{array}{c}0\\ K\\0\end{array}\right)$. This, on the other hand, means the parallel null plane $\cal H$ is given by
\[\cal H= \left\{ \left(\begin{array}{c}0\\X\\\s\end{array}\right)\mid X\in L, \s\in \rr\right\}.\]
But this was equivalent to properties \eqref{ray} and \eqref{null-ric}.
\eprf
\bbem
This tractorial characterisation in Theorem \ref{theo-null} yields integrability conditions in terms of the tractor curvature, which, on the other hand,  imply exactly the  obstructions \eqref{weyl} and \eqref{cotton} for the existence of an aligned pure radiation metric in a conformal class given in Theorem \ref{theo1}.  
Indeed, when splitting the tractor bundle with respect to the aligned pure radiation metric  $g_{ab}$  in the conformal class, the parallel null plane bundle $\cal H$ contains $K^B=\left(\begin{array}{c}0\\K^b\\
0\end{array}\right)$ and $\cal H^\bot$ contains $X^B=\left(\begin{array}{c}0\\X^b\\
0\end{array}\right)$ with $X^b$ orthogonal to $K^a$. Hence, with respect to a splitting in another metric $\hat g_{ab}=\e^{2\ups}g_{ab}$, using formula \eqref{splitting},  $K^B$ and $X^B$ are given as
\[K^B=
\left(
\begin{array}{c}
\ups_aK^a 
\\ K^b \\
0
  \end{array}\right),\ \ \
X^B=
\left(
\begin{array}{c}
\ups_aX^a 
\\ X^b   \\
0
  \end{array}\right)  ,
  \]
and still contained in $\cal H$ and $\cal H^\bot$, respectively. Hence, since $\cal H$
is parallel and thus invariant under the tractor curvature, applying the tractor curvature to $K^B$ 
yields a tractor that is again contained in $\cal H$. Pairing $\ol{R}_{abCD}K^C$ with $X^D$ and using \eqref{tracurv}, this implies that
\be
0&=&
\ol{R}_{abCD}K^CX^D\ =\ C_{abcd}K^cX^d,
\ee
 for all $X^d$ orthogonal to $K^c$, which is equation \eqref{weyl}.
  On the other hand, we have seen that $\cal H$ also contains the section $J^A=\left(\begin{array}{c}0\\0\\1\end{array}\right)$, when  splitting the tractor bundle with the aligned pure radiation metric $g_{ab}$. In another metric $\hat g_{ab}=\e^{2\ups}g_{ab}$, again using \eqref{splitting}, $J^A$ is given as
 \[J^A\ =\ 
\left(
\begin{array}{c}
- \einhalb  \ups^b\ups_b 
\\ \  \ups_a  \\
1 
  \end{array}\right).\]
  Since $\cal H^\bot$ is parallel and thus curvature invariant, pairing $\ol{R}_{abCD}X^C$ with $J^D$ yields
  \[
  0\ =\  
  \ol{R}_{abCD}X^CJ^D\ =\ 
  -A_{dab}X^d
+C_{ab\  c}^{\ \ d}X^c \ups_d,
\]
in which we have used formula \eqref{tracurv} for the tractor curvature. But this is
 equation \eqref{cotton} of Theorem \ref{theo1}.
  \ebem

Finally, from Theorem \ref{theo-null}, we get a characterisation of conformal pp--waves that is based on their characterisation as pure radiation metrics in Proposition \ref{ppprop}.

\bfolg
Let $(M,[g])$ be a pseudo-Riemannian conformal manifold of dimension $n>2$. 
Then, on an open and dense subset of $M$,  $[g]$ contains a pp--wave metric if and only if the normal conformal Cartan connection admits a parallel totally null plane subbundle $\cal H$ and the tractor curvature $\ol{R}$ satisfies
\begin{equation}
\ol{R}_{abCD}X^CY^D=0,\label{pptrac}
\end{equation}
for all $X^C$ and $Y^D$ orthogonal to $\cal H$.
\efolg
\bprf
If $g$ is a pp--wave in the conformal class, then by Proposition \ref{ppprop} it is an aligned  pure radiation metric, and hence, by Theorem \ref{theo-null} defines a parallel plane distribution in the tractor bundle. Using the pp--wave metric to split the tractor bundle, we have seen in the proof of the theorem that 
\[\cal H^\bot =
\left\{X^B=\left(\begin{array}{c}0\\X^b\\
\s\end{array}\right)\mid K_aX^a=0, \s\in \rr\right\}.\]
Hence, for $X^C=\left(\begin{array}{c}0\\X^c\\
\s\end{array}\right)$ and $Y^D=\left(\begin{array}{c}0\\Y^d\\
\tau\end{array}\right)$ from $\cal H^\bot$ with $X^c$ and $Y^d$ orthogonal to $K^a$, the 
 formula \eqref{tracurv} for the tractor curvature yields
 \begin{equation}\label{pptrac2}
\ol{R}_{abC}^{\ \ \ \ D}X^CY_D
\ =
\ -\tau A_{dab}X^d +\s A^c_{\ ab}Y_c +C_{ab\ \ d}^{\ \ c}X^dY_c
\ =\ 
0,\end{equation}
because of Propositions \ref{pureprop} and \ref{ppprop}.

On the other hand, by Theorem \ref{theo-null} the existence of the parallel plane distribution yields  an aligned pure radiation metric $g$ in the conformal class. With this metric as a gauge  the curvature condition \eqref{pptrac} then spells out as equation \eqref{pptrac2}. Now, since $g$ is an aligned pure radiation metric, by Proposition \ref{pureprop} we have that the first two terms in \eqref{pptrac2} vanish and we are left with
\[C_{ab\ \ d}^{\ \ \ c}X^dY_c=0\]
for all $X^d$ and $Y^c$ orthogonal to $K^a$. 
Then, from Proposition \ref{ppweylprop}  we know that $g$ is indeed a pp--wave.
\eprf







\def\cprime{$'$}

\end{document}